\documentclass[11pt]{amsart}
\usepackage{amsmath,amssymb,amsthm}
\usepackage{xcolor}
\usepackage{booktabs}
\usepackage{needspace}
\usepackage{microtype}
\usepackage[colorlinks=true,citecolor=blue!45!black,linkcolor=blue!45!black,urlcolor=blue!45!black]{hyperref}
\usepackage{mathtools}

\newtheorem{theorem}{Theorem}[section]
\newtheorem{lemma}[theorem]{Lemma}

\theoremstyle{remark}
\newtheorem{remark}[theorem]{Remark}

\title[Jordan schemes]{A construction of proper Jordan schemes of rank $5$}
\author{Sho Suda}
\address{Department of Mathematics, National Defense Academy of Japan, 1-10-20 Hashirimizu, Yokosuka, Kanagawa 239-8686, Japan}

\email{ssuda@nda.ac.jp}
\date{\today}
\hypersetup{pdftitle={Jordan schemes},pdfauthor={Sho Suda}}
\begin{document}

\begin{abstract}
Let $q$ be a prime power such that $m=q^2+q+1$ is a prime satisfying $m\equiv 1\pmod{4}$.  
We construct a proper symmetric Jordan scheme of rank $5$ on $2m$ vertices using a  Singer difference set with parameters $(m,q+1,1)$ and the Paley graph on $m$ vertices.  
\end{abstract}
\maketitle
\raggedbottom

\section{Introduction}
A symmetric Jordan scheme is a generalization of a symmetric association scheme, in which closure under ordinary matrix multiplication is replaced by closure under the Jordan product
\[
A\star B=\frac12(AB+BA).
\]
A Jordan scheme is called \emph{proper} if it is not the symmetrization of a coherent configuration. In particular, non-associativity does not imply properness.
Cameron \cite{Cam03} posed the following question: Does a proper Jordan scheme exist?

In \cite{KMR19}, Klin, Muzychuk, and Reichard
showed the existence of proper Jordan schemes.
Their examples include Jordan schemes on
$15$ and $24$ vertices of rank $5$ and on $40$ vertices of rank $7$. 
In \cite{MRK22}, Muzychuk, Reichard, and Klin
constructed several infinite families of proper
Jordan schemes and developed the basic theory
of Jordan schemes. 
They showed that the smallest possible rank of a proper symmetric Jordan scheme is $5$, and that its adjacency Jordan algebra is isomorphic to 
\[
\mathbb R\oplus\mathbb R\oplus\operatorname{Sym}_2(\mathbb R). 
\]
In \cite{HY25}, Hanaki and Yoshikawa constructed
two new proper Jordan schemes on $16$ and $32$
vertices related to quaternion and octonion
algebras. Their Jordan algebras have simple
components of non-Hermitian type.
In \cite{MPW26}, Muzychuk, Pech, and Woldar
classified thin Jordan schemes and established
a correspondence between regular thin Jordan
schemes and a special class of Moufang loops
called RA-loops.

Known examples of proper Jordan schemes of rank $5$ are as follows. 
Here, $(v;k_0,\ldots,k_4)$ denotes the number $v$ of vertices and the valencies $k_0,\ldots,k_4$. 
\begin{itemize}
    \item Small examples \cite{KMR19}:
    \[
    (15;1,2,4,4,4),\qquad
    (24;1,2,7,7,7).
    \]

    \item Switching construction \cite{KMR19,MRK22}:
    \[
    (3(r+1);1,2,r,r,r),
    \]
    where $r$ is a prime power satisfying
    $r\equiv1\pmod3$.

    \item Wallis--Fon-Der-Flaass construction \cite{MRK22}:
    \[
    \left(
    \frac{3^d(3^d+1)}{2};
    1,3^d-1,
    \frac{3^{d-1}(3^d-1)}{2},
    \frac{3^{d-1}(3^d-1)}{2},
    \frac{3^{d-1}(3^d-1)}{2}
    \right),
    \]
    where $d\ge2$ is even.
\end{itemize}

In this paper, we construct proper Jordan schemes
of rank $5$ on $2m=2(q^2+q+1)$ vertices,
where $q$ is a prime power and
$m=q^2+q+1$ is a prime satisfying
$m\equiv1\pmod4$.
Our construction uses a Singer difference set
with parameters $(m,q+1,1)$ and
the Paley graph on $m$ vertices.
The resulting Jordan schemes have the following parameters
    \[
    \left(
    2(q^2+q+1);
    1,\frac{q(q+1)}{2},
    \frac{q(q+1)}{2},q+1,q^2
    \right). 
    \]
These parameters differ from those of the rank~$5$ Jordan schemes
listed above.

For $q=3$, take the Singer difference set
$D=\{0,1,3,9\}\subseteq\mathbb Z_{13}$.
The union of two of our basic relations then gives the strongly
regular graph with parameters $(26,10,3,4)$ constructed by
Adel'son-Vel'skii, Weisfeiler, Leman and Faradzhev \cite{AWLF69}.
In Example~2.6 of \cite{dRJ92}, de Resmini and Jungnickel
describe this graph in terms of the Paley graph on $13$ vertices
and the difference set $D$.

\section{Preliminaries}
Let $X$ be a finite set and let $R_0,\ldots,R_s$ be nonempty subsets of $X\times X$. 
Let $A_i$ be the adjacency matrix of the directed graph $(X,R_i)$. 
The pair $(X,\{R_i\}_{i=0}^s)$ is said to be a symmetric \emph{Jordan scheme} if the following hold: 
\begin{enumerate}
    \item $A_0=I$;
    \item $\sum_{i=0}^s A_i=J$;
    \item $A_i^\top=A_i$ for $i=1,\ldots,s$; 
    \item for any $i,j,k$, 
    there exist coefficients $p_{i,j}^k\in\frac12\mathbb Z$ such that 
    \[
    A_i\star A_j:=\frac{1}{2}(A_iA_j+A_jA_i)=\sum_{k=0}^s p_{i,j}^k A_k. 
    \]
\end{enumerate}

Let $\mathcal J$ be the vector space spanned by
$A_0,\ldots,A_s$ over $\mathbb R$.

A subspace $\mathcal C$ of $M_X(\mathbb R)$
is called a coherent algebra if it contains
$I$ and $J$ and is closed under ordinary
matrix multiplication, transposition, and
the entry-wise product
\[
(A\circ B)_{xy}=A_{xy}B_{xy}.
\]
The symmetrization of $\mathcal C$ is defined by
\[
\operatorname{Sym}(\mathcal C)
=\{A\in\mathcal C:A^\top=A\}.
\]
A Jordan scheme is called \emph{proper} if
$\mathcal J\neq\operatorname{Sym}(\mathcal C)$
for every coherent algebra $\mathcal C$.

\section{Construction}
Let $q$ be a prime power and assume that
$m=q^2+q+1$ is a prime satisfying $m\equiv1\pmod4$.  

Let $D\subseteq\mathbb Z_m$ be a Singer difference set
with parameters $(m,q+1,1)$.
Let $U$ be the cyclic permutation matrix of order $m$,
and define
\[
N=\sum_{d\in D}U^d.
\]
Let $Q$ be the set of nonzero quadratic residues
in $\mathbb Z_m$, and define
\[
B=\sum_{a\in Q}U^a,\qquad C=J-I-B.
\]
Then $B$ is the adjacency matrix of the Paley graph
on $m$ vertices. 

Define $(0,1)$-matrices $A_0,\ldots,A_4$ by
\begin{align*}
    A_0&=\begin{bmatrix}
        I_m & O \\
        O & I_m
    \end{bmatrix},\quad
    A_1=\begin{bmatrix}
        B & O \\
        O & C
    \end{bmatrix},\quad
    A_2=\begin{bmatrix}
        C & O \\
        O & B
    \end{bmatrix},\\
    A_3&=\begin{bmatrix}
        O & N \\
        N^\top & O
    \end{bmatrix},\quad
    A_4=\begin{bmatrix}
        O & J-N \\
        J-N^\top & O
    \end{bmatrix}.
\end{align*}

Let $\mathfrak J$ be the configuration
with adjacency matrices $A_0,\ldots,A_4$.

\begin{theorem}
    $\mathfrak J$ is a symmetric Jordan scheme. 
\end{theorem}
\begin{proof}
    It is obvious that the conditions (1)-(3) of the Jordan scheme hold. We verify that the condition (4) holds. 
    Let $\mathcal J$ be the vector space spanned by $A_0,\ldots,A_4$ over the real number field $\mathbb{R}$.
    Since $B$ is the adjacency matrix of the Paley graph and $C=J-I-B$, 
    \begin{align*}
B^2 &=
\frac{q(q+1)}{2}I
+\frac{q(q+1)-4}{4}B
+\frac{q(q+1)}{4}C,\\
C^2 &=
\frac{q(q+1)}{2}I
+\frac{q(q+1)}{4}B
+\frac{q(q+1)-4}{4}C,\\
BC=CB &=
\frac{q(q+1)}{4}(B+C),\\
BJ=CJ &=
\frac{q(q+1)}{2}J.
\end{align*}
Then 
\begin{align*}
A_1\star A_1
&=\frac{q(q+1)}{2}A_0
+\frac{q(q+1)-4}{4}A_1
+\frac{q(q+1)}{4}A_2,\\
A_1\star A_2
&=A_2\star A_1=\frac{q(q+1)}{4}(A_1+A_2),\\
A_2\star A_2
&=\frac{q(q+1)}{2}A_0
+\frac{q(q+1)}{4}A_1
+\frac{q(q+1)-4}{4}A_2.
\end{align*}
Therefore $A_i \star A_j\in \mathcal  J$ for $i,j\in\{1,2\}$. 

Next, since $N$ is the incidence matrix of a Singer difference set,
\begin{align*}
NN^\top=N^\top N &= qI+J,\\
NJ=JN=N^\top J=JN^\top 
&= (q+1)J.
\end{align*}
Therefore 
\begin{align*}
A_3A_3
&=(q+1)A_0+A_1+A_2,\\
A_3A_4
&=A_4A_3=q(A_1+A_2),\\
A_4A_4
&=q^2A_0+q(q-1)(A_1+A_2),
\end{align*}
which implies that $A_i \star A_j\in \mathcal  J$ for $i,j\in\{3,4\}$. 

Since $B$ and $N$ are polynomials in the
    same cyclic permutation matrix $U$,
    they commute.
%    Indeed, letting
%    $B=\sum_{a\in Q}U^a$ and
%    $N=\sum_{d\in D}U^d$, we have
%    \[
%    BN=\sum_{a\in Q}\sum_{d\in D}U^{a+d}
%      =\sum_{d\in D}\sum_{a\in Q}U^{d+a}
%      =NB.
%    \]
    Thus $N,N^\top,B$, and $C=J-I-B$
    commute pairwise, and
    \begin{align*}
    BN+NC
    &=(B+C)N
     =(J-I)N
     =(q+1)J-N,\\
    CN^\top+N^\top B
    &=(B+C)N^\top
     =(J-I)N^\top
     =(q+1)J-N^\top.
    \end{align*}
    Therefore,
    \begin{align*}
    A_1\star A_3=A_2\star A_3
    &=\frac q2 A_3+\frac{q+1}{2}A_4,\\
    A_1\star A_4=A_2\star A_4
    &=\frac{q^2}{2}A_3
      +\frac{q^2-1}{2}A_4.
    \end{align*}
    Hence $A_i\star A_j\in\mathcal J$
    for $i\in\{1,2\}$ and $j\in\{3,4\}$.

    This verifies that condition (4) holds,
    and thus the set $\{A_0,\ldots,A_4\}$
    is a Jordan scheme.
\end{proof}

\begin{remark}\label{rem:1}
    We show that the Jordan scheme is not an
association scheme.
Indeed,
\[
A_1A_3-A_3A_1
=
\begin{bmatrix}
O & N(B-C)\\
-(N(B-C))^\top & O
\end{bmatrix}.
\]
Since $NN^\top=qI+J$, the matrix $N$ is
nonsingular. Moreover, $B\ne C$.
Hence $A_1A_3\ne A_3A_1$.
Since all $A_i$ are symmetric, their
linear span is not closed under ordinary
matrix multiplication.
Thus the Jordan scheme is not
an association scheme.
\end{remark}

In the $26$-vertex example described in the introduction,
the strongly regular graph has adjacency matrix $A_1+A_3$.
The two copies of $\mathbb Z_{13}$ induce the Paley graph and
its complement, giving the exceptional strongly regular
decomposition described by Haemers and Higman \cite[p.~393]{HH89}.
For this decomposition, the five adjacency matrices in
\cite[Eq.~(7.6)]{vDMM13} are, in the order given there,
$A_0,A_3,A_1,A_4,A_2$.
Since the decomposition is exceptional,
Proposition~7.7 of \cite{vDMM13} also shows that these relations
do not form an association scheme.

\section{Properness}
In order to prove the properness for the Jordan scheme, we use the following. 
\begin{lemma}\label{lem:1}
    Let $K=A_iA_j-A_jA_i$. If $K\circ K$ is not in $\mathcal J$, then the Jordan scheme $\mathfrak J$ is proper.
\end{lemma}
\begin{proof}
    Assume, to the contrary, that $\mathfrak J$ is not proper. 
    By definition, the Jordan algebra $\mathcal J$ is equal to $\operatorname{Sym}(\mathcal C)$, where $\mathcal C$ is the coherent algebra of a coherent configuration. 

    Since $K\in \mathcal C$ and $\mathcal C$ is closed under the entry-wise product, $K\circ K\in \mathcal C$. 
    Since $K$ is skew-symmetric, $K \circ K$ is symmetric. 
    So $K\circ K\in \operatorname{Sym}(\mathcal C)=\mathcal J$, a contradiction to the assumption that $K\circ K$ is not $\mathcal J$. 
\end{proof}
\begin{theorem}
    The Jordan scheme $\mathfrak J$ is proper. 
\end{theorem}
\begin{proof}
    Let $K=A_1A_3-A_3A_1$. We verify that the matrix $K$ satisfies the assumption of Lemma~\ref{lem:1}. 
    Let $S(t)$ satisfy $N(B-C)=\sum_{t=0}^{m-1}S(t)U^t$. Then 
    \begin{align*}
        S(t)=\sum_{d\in D}\chi(d-t),
    \end{align*}
    where $\chi$ is the quadratic character. 
    Since $D$ is a difference set, 
    \[
    \sum_{t\in D}S(t)=\sum_{t,d\in D}\chi(d-t)=\sum_{x=0}^{m-1}\chi(x)=0.
    \]
    With 
    \[
    \sum_{t=0}^{m-1}S(t)=\sum_{d\in D}\sum_{t=0}^{m-1}\chi(d-t)=\sum_{d\in D}\sum_{x=0}^{m-1}\chi(x)=0, 
    \]
    we obtain 
    \begin{align*}\label{eq:1}        
    \sum_{t\not\in D}S(t)=0.
    \end{align*}

%    We show that either $i=3,4$, there exist pairs of vertices $(x,y),(x',y')\in R_i$ such $K_{x,y}\neq K_{x',y'}$, that is, there exist integers $t,t'\in D$ such $|S(t)|\neq |S(t')|$, or $t,t'\not\in D$ such $|S(t)|\neq |S(t')|$.  
    We show that $|S(t)|$ is not constant on at least one of $D$ and its complement, which implies that $K\circ K\not\in \mathcal J$. 
    Assume, to the contrary, that $|S(t)|$ depends only on whether $t\in D$ or not. Let $\alpha=|S(t)|$ for $t\in D$ and $\beta=|S(t)|$ for $t\not\in D$.   

    (i): 
Suppose that $q$ is even.
Since $|D|=q+1$ is odd and
$\sum_{t\in D}S(t)=0$, we have $\alpha=0$.
The orthogonality of the quadratic character gives
\[
q^2\beta^2=q^2(q+1),
\]
and hence $\beta^2=q+1$.
Since $q$ is a power of $2$, $(\beta+1)(\beta-1)=q=2^f$ for some positive integer $f$. 
Since $\beta-1$ and $\beta+1$ are powers of $2$ differing by $2$, we obtain  
$\beta=3$ and $q=8$.

In this case, $m=73$. 
Since $2$ is a multiplier of the Singer difference set \cite[Corollary 1.2]{JV84},
we may translate $D$ so that $2D=D$.
Since $2$ has multiplicative order $9$ modulo $73$,
every nonzero orbit of $\langle 2\rangle$ has size $9$.
As $|D|=9$, we have
$D=a\langle 2\rangle$
for some $a\in\mathbb Z_{73}^\times$.

Since $\langle 2\rangle\subseteq Q$,
we have
\[
|S(0)|=9,
\]
contrary to $\beta=3$.

    (ii): 
Suppose that $q$ is odd.
Since $|\mathbb Z_m\setminus D|=q^2$ is odd and
$\sum_{t\notin D}S(t)=0$, we have $\beta=0$.
 
By the orthogonality of the quadratic character,
\[
\sum_{t=0}^{m-1}S(t)^2=q^2(q+1).
\]
Then $(q+1)\alpha^2=q^2(q+1)$, which implies that $\alpha=q$.
Since $\sum_{t\in D}S(t)=0$, there exist
$x,y\in D$ such that $S(x)=q$ and $S(y)=-q$.

Since
\[
S(x)=\sum_{d\in D\setminus\{x\}}\chi(d-x)=q,
\]
we have $\chi(d-x)=1$ for any
$d\in D\setminus\{x\}$.
Similarly, $\chi(d-y)=-1$ for any
$d\in D\setminus\{y\}$.
In particular,
\[
\chi(y-x)=1,\qquad \chi(x-y)=-1,
\]
which contradicts $\chi(-1)=1$ since
$m\equiv1\pmod4$.
  
    Therefore, $K\circ K$ is not in $\mathcal J$. By Lemma~\ref{lem:1}, the theorem is proved. 
\end{proof}

\section*{Acknowledgments}
This work was supported by JSPS KAKENHI Grant Number JP26K06904.

\section*{Declaration of AI use}
The author used ChatGPT (OpenAI) to assist with English language
editing, checking mathematical arguments, identifying relevant
references, and comparing earlier results with those of this paper.
The author reviewed the output and takes full responsibility for
the paper and its references.

\end{document}